\documentclass{siamart190516}

\usepackage[utf8]{inputenc}
\usepackage[english]{babel}

\usepackage{amsmath,amsfonts,amssymb,bm}

\usepackage{algpseudocode}
\usepackage{multirow}
\usepackage[per-mode = symbol]{siunitx}
\DeclareSIUnit{\flop}{flop}
\DeclareSIUnit{\byte}{B}

\usepackage[caption=false]{subfig}
\usepackage{graphicx}
\usepackage{array}
\usepackage[caption=false]{subfig}

\usepackage[textsize=footnotesize]{todonotes}

\usepackage{mdframed}
\usepackage{booktabs}
\usepackage{xspace}

\newmdenv[backgroundcolor=lightgray,frametitle=Todo]{todobox}
\DeclareMathOperator {\vspan}{span}

\usepackage{xstring}
\newcommand{\formatindex}[1]{\StrSubstitute{#1}{:}{\,\mathtt{:}\,}}
\newcommand{\vslice}[1]{\texttt{[}\formatindex{#1}\texttt{]}}
\newcommand{\mslice}[2]{\texttt{[}\formatindex{#1},\formatindex{#2}\texttt{]}}
\newcommand{\cpp}{C\texttt{++}\xspace}
\newcommand{\methodname}[1]{\texttt{#1}}

\newcommand{\thetitle}{A Hardware-aware and Stable Orthogonalization Framework}

\title{\thetitle
  \thanks{Submitted to the editors April 28, 2021.
\funding{Funded by the Deutsche Forschungsgemeinschaft (DFG, German Research Foundation) under Germany’s Excellence Strategy EXC 2044-390685587, Mathematics
Münster: Dynamics-Geometry-Structure}}}
\author{Nils-Arne Dreier\thanks{University of Münster,
    Orléansring 10, 48149 Münster (\email{n.dreier@uni-muenster.de}, \email{christian.engwer@uni-muenster.de})}
  \and Christian Engwer\footnotemark[2] }

\headers{\thetitle}{N.-A. Dreier, C. Engwer}

\ifpdf
\hypersetup{
  pdftitle={\thetitle},
  pdfauthor={N.-A. Dreier, C. Engwer}
}
\fi

\begin{document}
\maketitle
\begin{abstract}
  The orthogonalization process is an essential building block in Krylov
  space methods, which takes up a large portion of the computational time.
  Commonly used methods, like the Gram-Schmidt method,
  consider the projection and normalization separately and store the orthogonal
  base explicitly.
  We consider the problem of orthogonalization and normalization
  as a QR decomposition problem on which we apply known algorithms,
  namely CholeskyQR and TSQR.
  This leads to methods that solve the orthogonlization problem with reduced
  communication costs, while maintaining stability and stores the orthogonal
  base in a locally orthogonal representation.
  Furthermore, we discuss the novel method as a framework which allows us to
  combine different orthogonalization algorithms and use the best algorithm
  for each part of the hardware.
  After the formulation of the methods, we show their advantageous performance
  properties based on a performance model that takes data transfers within
  compute nodes as well as message passing between compute nodes into account.
  The theoretic results are validated by numerical experiments.
\end{abstract}

\begin{keywords}
  Orthogonalization, Block Krylov methods, High-Performance Computing
\end{keywords}

\begin{AMS}
  15A23, 65F25, 65Y05
\end{AMS}

\section{Introduction}
The orthogonalization process is an important building block in
Krylov space methods, both to solve linear systems as well as to compute eigenvectors.
In this paper we focus on the orthogonalization as part of the Arnoldi process,
as it is used for example in the GMRes method.
Beside the application of the operator (and the preconditioner) the projection
and orthogonalization step takes up the major portion of the runtime.
On modern architectures not the actual computation is the bottleneck, but
communication~\cite{dongarra2011international}.
Communication means the exchange of data between the components of the hardware.
This exchange happens between the memory and the CPU as well as between compute
nodes in distributed environments.
This paper aims at algorithms that optimize the communication and hence achieve
better performance in high-performance computing (HPC) environments.
To this end, we extend the TSQR algorithm
\cite{demmel2008communication,demmel2012communication} and adopt the data
structure to store the base of the Krylov space in a locally orthogonal representation.
Based on it we introduce a
framework in which different algorithms can be combined utilizing their
advantages on different parts of the computer architecture.

Communication-avoiding and communication-hiding Krylov methods gained a lot attention
these days.
In particular Block Krylov methods are well suited for high-performance computing
\cite{dreier2019strategies, dreier2020thesis}.
These methods were originally developed to solve linear systems with multiple
right-hand sides \cite{oleary1980block} or to compute multiple eigenvectors
simultaneously \cite{golub1977block}.
Since then also block variants of the Arnoldi method \cite{sadkane1993block} and
GMRes method \cite{simoncini1996convergence} were proposed.

A particular issue of the orthogonalization process is stability.
Especially the classical Gram-Schmidt procedure may produce a significant
orthogonalization error.
The stability properties of various block Gram-Schmidt procedures are analyzed
in the paper of Carson et al.~\cite{carson2020stability} and an excellent
overview of existing block Gram-Schmidt methods is given in
\cite{carson2020overview}, taking performance and stability into account.

A new approach for reducing the synchronization cost in a block Gram-Schmidt
algorithm are presented by Swirydowicz et al.~\cite{swirydowicz2021low}.
To remedy the problem of synchronization points in the block Gram-Schmidt
process, they propose the so-called
low-synch~\cite{swirydowicz2021low,bielich2021lowsynch}
methods that delay the normalization of Krylov vectors.

Yamazaki et al.~\cite{yamazaki2019low} already presented a method that combines
the projection and normalization step in a single reduction in the context of
the classical Gram-Schmidt method.
This was later extended to block Gram-Schmidt \cite{carson2020overview} under the name \methodname{BCGS-PIP}.
The same algorithm can be derived in our setting as we will see in \cref{sec:bcgs-pip}.

The separate consideration of projection and normalization leads to multiple
synchronization points, at least two.
This paper is based on the observation, that these problems can be treated
together:
Let $Q \in \mathbb{R}^{n\times k}$ be an
orthogonal matrix which columns span the already computed Krylov space and
$X\in\mathbb{R}^{n\times s}$ the block-vector that contains the new directions
of the Krylov space. The aim is now to compute an orthogonal matrix $U \in
\mathbb{R}^{n\times t}$ such that $\vspan(Q, U) = \vspan(Q, X)$ and $Q^TU = 0$.
We call this problem the \emph{project and normalize} (PQR) problem.
It can be solved by computing the reduced QR factorization
\begin{align}
  \label{eq:pqr}
  \begin{bmatrix}
    Q & X
  \end{bmatrix}
  &=
    \begin{bmatrix}
      Q & U
    \end{bmatrix}
          \begin{bmatrix}
            I & P\\
            0 & N
          \end{bmatrix},
\end{align}
with projection matrix $P\in\mathbb{R}^{k\times s}$ and normalizer
$N\in\mathbb{R}^{t\times s}$.
In most cases we have $t=s$ and $N$ is a upper triangular matrix.
The other case, $t<s$, appears if the system $\begin{bmatrix}Q&X\end{bmatrix}$
is linearly dependent.
In the context of block Krylov space methods this case is called the deflation case.
The matrix $N$ is then a row-echelon matrix.
Depending on the implementation details, it might be reasonable to use $t=s$ in
all cases and treat the rank-deficiency of $N$ later on. One possibility is to
use a rank-revealing QR-decomposition, i.e.\ pivoting, to enforce triangular shape of $N$.

The matrix $Q$ must not necessarily stored explicitly, but, depending on the
orthogonlization method, it might be useful to store it in an other representation.
Walker \cite{walker1988implementation} for example used a sequence of Householder
reflectors to store the orthogonal base in the GMRes method. This was proposed
by Walker \cite{walker1988implementation}.
In \cref{sec:tspqr} we will use a locally orthogonal representation.
However, the output $U$ of the method should be available explicitly, as it used
to compute the subsequent Krylov space directions by applying the operator.

Computing the QR-factorization \eqref{eq:pqr} is actually more costly than the
original orthogonalization problem,
but as we know QR-factorization algorithms (namely CholQR~\cite[Thm. 5.2.3]{golub2013matrix} and TSQR~\cite{demmel2012communication}) that can
solve this problem with only one global synchronization, we can apply these
algorithm to solve the projection and normalize problem with one global
synchronization.
Additionally, we can make use of the known structure of problem \eqref{eq:pqr},
i.e.\ $Q$ is already orthogonal, to simplify and improve the resulting method.

The main contribution of this paper is an orthogonalization
framework for Krylov methods so, that the actual algorithm can be adapted to the
properties of the hardware.
It is structured as follows.
In \cref{sec:sota} we show how a method that solves problem \eqref{eq:pqr} can be used
in the Arnoldi procedure and review the
Gram-Schmidt and Householder method for the orthogonalization.
These methods are the fundamental building blocks in the framework.
As a further building block we discuss the \methodname{BCGS-PIP} method,
and show how to deduce it
from the Cholesky QR algorithm in \cref{sec:bcgs-pip}.
In the same manner we introduce the novel \methodname{TreeTSPQR} and
\methodname{FlatTSPQR} methods, which are deduced from
the respective TSQR algorithm in \cref{sec:tspqr}, which act as connectors between the
building blocks.
In \cref{sec:framework} we use performance models to analyze the
performance of the different algorithms and show how
the TSPQR algorithms can be combined with the Householder and \methodname{BCGS-PIP+} method to build an
orthogonalization framework that is adapted on the architecture of a given
supercomputer.
Numerical experiments that show the stability and performance properties are
given in \cref{sec:numeric}.
Finally, we give a conclusion and outlook in \cref{sec:conclusions}.

\section{Orthogonalization in Krylov space methods}
\label{sec:sota}
Let $R\in\mathbb{R}^{n\times s}$ with ${s \ll n}$. The block Arnoldi method constructs a sequence of orthogonal bases $\mathcal{V}_k$
for the block Krylov space
\begin{align}
  \mathcal{K}^k(A, R) &= \vspan\left(R, AR, \ldots, A^{k-1}R\right)
\end{align}
and block-Hessenberg matrices $H_k$ that satisfy the so-called block Arnoldi
relation
\begin{align}
  A\mathcal{V}_{k-1} = \mathcal{V}_{k}H_k.
\end{align}

The algorithm can be formulated by solving problem \eqref{eq:pqr} in every
iteration.
The algorithm is shown in \cref{alg:blockarnoldi}.

\begin{algorithm}
  \caption{Block Arnoldi}
  \label{alg:blockarnoldi}
  \begin{algorithmic}
    \State $V_0 = RH_0$\Comment QR-factorization
    \For{$i=0,\ldots$}
    \State $X_{i+1} = AV_i$
    \State $\begin{bmatrix}\mathcal{V}_i& X_{i+1}\end{bmatrix}
    = \begin{bmatrix}\mathcal{V}_i & V_{i+1}\end{bmatrix}
    \begin{bmatrix}
      I & P\\
      0 & N
    \end{bmatrix}
    $ \Comment Solve \eqref{eq:pqr}
    \State $\mathcal{V}_{i+1} = \begin{bmatrix}
      \mathcal{V}_i & V_{i+1}
    \end{bmatrix}$
    \State $H_{i+1} =
    \begin{bmatrix}
      H_i & P\\
      0 & N
    \end{bmatrix}$
    \EndFor
  \end{algorithmic}
\end{algorithm}

The basis $\mathcal{V}_k$ can then be used in the block GMRes method to solve
large sparse linear systems with multiple right-hand sides or to compute
eigenvectors of the operator $A$.
The spectrum of $H_k$ is an approximation for the spectrum of $A$.

In practice two strategies are commonly used to solve \eqref{eq:pqr} in \cref{alg:blockarnoldi}.
The Gram-Schmidt method, either in the classical or in the modified
variant, and the Householder method.
Depending on the used orthogonalization method the matrices $\mathcal{V}_k$
might not be stored explicitly.
As we will see, if the Householder method is used, the basis $\mathcal{V}_k$ is stored as
product of Householder reflectors.
Similarly, in the method introduced in \cref{sec:tspqr} it is stored as the
product of locally orthogonal matrices.

To assess the stability of the presented algorithms we apply them for computing
a QR-factorization of a matrix $A$.
\cref{alg:qr} computes the QR-factorization by
solving the PQR problem \eqref{eq:pqr} block-column-wise, where $A_i$ denote the
$i$th block-column of $A$ i.e.\ the columns $si$ to $s(i+1)$.
\begin{algorithm}
  \caption{Block Column-oriented QR-factorization}
  \label{alg:qr}
  \begin{algorithmic}
    \State $A_0 = Q_0R_0$ \Comment{QR-factorization}
    \For{$i=1\,\ldots\,\frac{m}{s}$}
    \State $\begin{bmatrix} Q_{i} & A_{i+1}
    \end{bmatrix} = \begin{bmatrix}
      Q_{i} & U_u
    \end{bmatrix}
    \begin{bmatrix}
      I & P\\
      0 & N
    \end{bmatrix}$
    \Comment{Solve \eqref{eq:pqr}}
    \State $Q_{i+1} = \begin{bmatrix}
      Q_i & U_i
    \end{bmatrix}$
    \State $R_i = \begin{bmatrix}
      R_{i-1} & P\\
      0 & N
    \end{bmatrix}$
    \EndFor
  \end{algorithmic}
\end{algorithm}

The orthogonality error in the
Frobenius norm
\begin{align}
  \label{eq:ortho_error}
  e_\perp := \|I - Q^TQ\|_F
\end{align}
gives an indication for the stability, as accumulated errors immediately impede
the orthogonality and increase the error $e_\perp$.
Estimations for $e_\perp$ are given depended on the machine precision and condition number of the input matrix $A$,
cf.\ \cref{sec:numeric}.
We will denote the machine precision by $\varepsilon$ and the condition number
by $\kappa$.

\subsection{Gram-Schmidt}
The Gram-Schmidt method is the oldest variant for orthogonalizing a set for
vectors.
The block classical variant (\methodname{BCGS}) projects the matrix $X$
onto the orthogonal complement of $Q$ by computing
\begin{align}
  \tilde{U} = (I - QQ^T)X = X - Q(Q^TX).
\end{align}
After that the result is normalized by computing the reduced QR decomposition
of $\tilde{U}$.
This QR factorization can for example be efficiently computed by the TSQR
algorithm or by the CholQR algorithm.
The resulting algorithm is shown in \cref{alg:cgs}.
\begin{algorithm}
  \caption{Block Classical Gram-Schmidt (\methodname{BCGS})}
  \label{alg:cgs}
  \begin{algorithmic}
    \State $P = Q^TX$
    \State $\tilde{U} = X - QP$
    \State $\tilde{U} = UN$ \Comment QR factorization
  \end{algorithmic}
\end{algorithm}
As shown by Giraud et al.~\cite{giraud2005rounding} the orthogonalization
error of the \methodname{BCGS} method is of order
$\mathcal{O}(\varepsilon\kappa^2)$.

To improve the stability the modified Gram-Schmidt method was
introduced and adapted for block orthogonalization \cite{jalby1991stability}.
It computes the projection of each column $q_i$ in $Q$ separately
\begin{align}
  \tilde{U} = (I-q_kq_k^T)\cdots(I-q_1q_1^T)X.
\end{align}
This leads to an orthogonalization error of $\mathcal{O}(\varepsilon\kappa)$.
\cref{alg:mgs} shows the pseudocode of the algorithm. Subscripts denote the
column of the respective matrix.
\begin{algorithm}
  \caption{Block Modified Gram-Schmidt (\methodname{BMGS})}
  \label{alg:mgs}
  \begin{algorithmic}
    \State $\tilde{U}_0 = X$
    \For{$i=1,\ldots,k$}
    \State $P_i = q_i^T\tilde{U}_{i-1}$
    \State $\tilde{U}_i = \tilde{U}_{i-1} - q_iP_i$
    \EndFor
    \State $\tilde{U}_k = UN$\Comment QR factorization
  \end{algorithmic}
\end{algorithm}

The modified Gram-Schmidt method has the disadvantage that $k$ inner products
are computed sequentially.
In contrast to the classical variant this means that the matrix $U$ is loaded
$k$ times from the memory and in a distributed computation environment this
leads to $k$ collective communications, i.e.\ synchronization points.

Another approach for improving the stability is the reiteration of the classical
Gram-Schmidt method \cite{abdelmalek1971round}.
It improves the orthogonalization error of the classical Gram-Schmidt method by
applying the algorithm twice.
It is proven that the first reiteration brings the orthogonalization error down to
$\mathcal{O}(\varepsilon)$ \cite{barlow2013reorthogonalized}.
\cref{alg:cgsi} shows the reiterated \methodname{BCGS} algorithm (\methodname{BCGS+}).
As the \methodname{BCGS} method is applied twice, the costs for the \methodname{BCGS+} method are twice
the costs for the \methodname{BCGS} method.

\begin{algorithm}
  \caption{Block Classical Gram-Schmidt with reiteration (\methodname{BCGS+})}
  \label{alg:cgsi}
  \begin{algorithmic}
    \State Solve $
    \begin{bmatrix}
      Q&X
    \end{bmatrix} =
    \begin{bmatrix}
      Q&U_1
    \end{bmatrix}
    \begin{bmatrix}
      I & P_1\\
      0 & N_1
    \end{bmatrix}$
    using \methodname{BCGS}
    \State Solve $
    \begin{bmatrix}
      Q&U_1
    \end{bmatrix} =
    \begin{bmatrix}
      Q&U
    \end{bmatrix}
    \begin{bmatrix}
      I & P_2\\
      0 & N_2
    \end{bmatrix}$
    using \methodname{BCGS}
    \State $P = P_1N_2 + P_2$ and $N = N_2 + N_1N_2$
  \end{algorithmic}
\end{algorithm}

\subsection{Householder}
We consider the Householder method in a form that treats the problem column-wise,
an elaborate analysis can be found in the book of Golub and von Loan
\cite[Chapter 5]{golub2013matrix}.
For brevity we employ a python-style slicing syntax to refer to submatrices
and subvectors. It is indicated by rectangular brackets. E.g.\ $X\mslice{0:k}{i}$ denotes
the first $k$ entries of the $i$th column of $X$.

The matrix $Q$ in the Householder method is stored as a product of Householder
reflectors
\begin{align}
  \label{eq:hh}
  \mathcal{Q}_k &= H_0\cdots H_{k-1}&\text{with}&&H_i = I - 2v_iv_i^T,
\end{align}
for some $v_i\in\mathbb{R}^n$ with $\|v_i\| = 1$.
The operators $H_i$ are symmetric and orthogonal.
The matrix $Q$ can be explicitly assembled by applying the operator
$\mathcal{Q}_k$ to the first $k$ unit vectors
\begin{align}
  Q = \mathcal{Q}_kI\mslice{0:n}{0:k}.
\end{align}

To solve problem \eqref{eq:pqr}, the new Householder reflectors are constructed
column by column. Let $\mathcal{Q}_{k+i} = H_0\cdots H_{k+i-1}$ be
constructed.
The next column $X\mslice{0:n}{i}$ in
$X$ is projected by applying $\mathcal{Q}^T_{k+i}$
\begin{align}
  \tilde{u}_i = \mathcal{Q}_{k+i}^TX\mslice{0:n}{i} = H_{k+i-1}\cdots H_0 X\mslice{0:n}{i}.
\end{align}
The upper $k$ entries in $\tilde{u}_i$ are the respective coefficients of $P$ and
$N$
\begin{align}
  P\mslice{0:n}{i} &= \tilde{u}_i\vslice{0:k} & N\mslice{0:i}{i} = \tilde{u}_i\vslice{k:k\!+\!i}.
\end{align}
From the lower part the vector $\tilde{v}_{k}$ is computed as
\begin{align}
  \tilde{v}_{k} = (\tilde{u}\vslice{k\!+\!i:n} - \sigma\lambda e_1),
\end{align}
where $\lambda$ is the norm of the lower part of $\tilde{u}$ and $\sigma$ chosen
as the invert sign of $\tilde{u}\vslice{k}$ to avoid cancellation and improve the
stability of the algorithm.
The diagonal entry of $N$ is then given by $N\mslice{i}{i} = \sigma\lambda$.
This vector is then normalized to obtain the Householder reflector
\begin{align}
  v_{k+1}\vslice{k:n} = \frac{1}{\|\tilde{v}_{k+1}\|}\tilde{v}_{k+1},
\end{align}
while the upper $k$ coefficients are set to $0$.

\begin{algorithm}
  \caption{Householder (\methodname{HH})}
  \label{alg:householder}
  \begin{algorithmic}[1]
    \State $\tilde{X} = H_{k-1}\cdots H_0 X$
    \For{$i=0\ldots s-1$}
    \State $\tilde{u}_i = H_{k+i-1} \cdots H_k \tilde{X}\mslice{0:n}{i}$
    \State $P\mslice{0:n}{i} = \tilde{u}_i\mslice{0:k}{i}$
    \State $N\mslice{0:i}{i} = \tilde{u}_i\mslice{k:k\!+\!i}{i}$
    \State $\lambda = \|\tilde{u}_i\mslice{k:n}{i}\|_2$
    \State $\sigma = -\operatorname{sgn}(\tilde{u}_i\mslice{k}{i})$
    \State $N\mslice{i}{i} = \sigma\lambda$
    \State $\tilde{v}_{k+i} = \tilde{u}\mslice{k\!+\!i:n}{i} - \sigma\lambda e_1$
    \State $v_{k+i}\vslice{k\!+\!i:n} = \frac{1}{\|\tilde{v}_{k+i}\|}\tilde{v}_{k+i}$
    \EndFor
  \end{algorithmic}
\end{algorithm}
\cref{alg:householder} shows the pseudocode for the Householder
method (\methodname{HH}).
In practice the application of the Householder reflectors $H_{k-1}\cdots H_0$ can be applied on
block for all columns in $X$.
This avoids reading $v_k\ldots v_1$ from the memory in every loop iteration.
The algorithm is known to be stable and produces and orthogonalization error in $\mathcal{O}(\varepsilon)$.

\section{BCGS-PIP}
\label{sec:bcgs-pip}
The CholQR algorithm is a well-known algorithm for computing the QR
factorization of a tall-skinny matrix.
It computes the QR factorization of a matrix $A$ by Cholesky factorizing the
Gram matrix $A^TA = LL^T$.
The QR factorization is then given by $A = QR$ with $Q = XL^{-T}$ and $R=L^T$.
While this algorithm has optimal performance properties, it renders 
unstable if the matrix $A$ is ill-conditioned.
As it was shown by Carson et al.~\cite{carson2020stability} and we will also see in \cref{sec:numeric} the orthogonalization error is of
order $\mathcal{O}(\kappa^2 \varepsilon)$, as long as $\kappa^2 < \varepsilon^{-1}$.

To deduce the algorithm, we apply the Cholesky QR algorithm on
the PQR problem \eqref{eq:pqr}. The occurring Gram matrix can be Cholesky factorized by
\begin{align}
  \begin{bmatrix}
    Q & X
  \end{bmatrix}^T
        \begin{bmatrix}
          Q & X
        \end{bmatrix} &=
                        \begin{bmatrix}
                          I & P\\
                          P^T & X^TX
                        \end{bmatrix}\\
      &=
        \begin{bmatrix}
          I & 0\\
          P^T & N^T
        \end{bmatrix}
                \begin{bmatrix}
                  I & P\\
                  0 & N
                \end{bmatrix},
\end{align}
where $P=Q^TX$ and $N$ is the Cholesky factor of $X^TX - P^TP = N^TN$.
The result $U$ is then given by
\begin{align}
  U = XL^{-T} - QPL^{-T}.
\end{align}
This algorithm is called \methodname{BCGS-PIP} \cite{carson2020stability}. The PIP stands for
Pythagorean Inner Product and refers to the original derivation from the
Pythagorean theorem.
As it already holds for the CholQR algorithm, this algorithm turns out to be
quite unstable.
Like CholQR or \methodname{BCGS}, the orthogonlization error is of order
$\mathcal{O}(\kappa^2 \varepsilon)$.
For a detailed stability analysis see \cite{carson2020overview}.
The pseudocode for this algorithm is shown in Algorithm \ref{alg:bcgs-pip}.

\begin{algorithm}
  \caption{\methodname{BCGS-PIP}}
  \label{alg:bcgs-pip}
  \begin{algorithmic}
    \State Compute $P = Q^TX$ and $G = X^TX$ \Comment{At once}
    \State Compute Cholesky factorization $N^TN = G - P^TP$
    \State $U = XN^{-1} - QPN^{-1}$
  \end{algorithmic}
\end{algorithm}

Regarding the communication, the algorithm requires one synchronization for
computing the Gram and projection matrix.
The matrices $Q$ and $X$ are loaded two times from the memory, once for the
computation of $G$ and $P$ and once for the assembly of $U$.

Using this algorithm in the GMRes method leads to the so-called \textit{one-step
  latency method} presented by Ghysels et al.~\cite{ghysels2013hiding}.

\subsection{Reiteration}
A remedy for the loss of orthogonalization is to
reiterate the algorithm \cite{yamamoto2015roundoff}.
As for the \methodname{BCGS+} method, this brings the orthogonalization
error down to machine precision as long as $\varepsilon\kappa^2 \leq \frac12$,
while doubling the computational and communication effort.
In total this algorithm performs two synchronizations and loads the matrices $Q$
and $X$ four times.
\methodname{BCGS-PIP} with reiteration (\methodname{BCGS-PIP+}) is
shown in \cref{alg:bcgsi-pip}.
\begin{algorithm}
  \caption{\methodname{BCGS-PIP} with reiteration (\methodname{BCGS-PIP+})}
  \label{alg:bcgsi-pip}
  \begin{algorithmic}
    \State Compute $P_1 = Q^TX$ and $G_1 = X^TX$ \Comment{At once}
    \State Compute Cholesky factorization $N_1^TN_1 = G_1 - P_1^TP_1$
    \State $U_1 = XN_1^{-1} - QP_1N_1^{-1}$
    \State Compute $P_2 = Q^TU_1$ and $G_2 = U_1^TU_1$ \Comment{At once}
    \State Compute Cholesky factorization $N_2^{T}N_2 = G_2 - P_2^TP_2$
    \State $U = U_1N_2^{-1} - QP_2N_2^{-1}$
    \State $P = P_1N_2 + P_2$ and $N = N_2 + N_1N_2$
  \end{algorithmic}
\end{algorithm}

Another drawback of the algorithm is the implementation of deflation.
The textbook variant of the Cholesky factorization is not rank-revealing,
meaning it can not be used to decide whether $G$ has full rank.
Hence we perform a singular value-decomposition of the matrix $G-P^TP = U \Sigma
U^T$, where we truncate singular values smaller than a certain tolerance,
yielding a full-rank square-root
\begin{align}
  N = \tilde{U}\tilde{\Sigma}^{\frac12}
\end{align}
of $G-P^TP$, where $\tilde{\Sigma}$ denotes the diagonal matrix of the
significant singular values and $\tilde{U}$ the corresponded left columns of $U$.

\section{TSPQR}
\label{sec:tspqr}
After we have introduced several orthogonalization methods, we
now propose the TSPQR method.
It subdivides the PQR problem into smaller problems of the same type and thus
enables us to apply recursion.
Furthermore, it builds up on the previous introduced methods and gives us the
choice which method is used on which level.

Like the \methodname{BCGS-PIP} algorithm the TSPQR algorithm is derived from the
QR-factorization problem \eqref{eq:pqr}, but instead of the CholQR algorithm, we
apply the TSQR algorithm \cite{demmel2012communication}.
It was proposed in two variants - a tree reduction and a sequential incremental
one.
We will consider both in the following subsections.
\subsection{TreeTSPQR}
The tree-TSQR algorithm computes a QR factorization of a tall-skinny matrix $A\in
\mathbb{R}^{n\times s}$ by row-wise decomposing $A$ into $p$ \emph{local} blocks and
factorizing these blocks in parallel
\begin{align}
  \label{eq:localqr}
  A = \begin{bmatrix}
    A_1\\
    \vdots\\
    A_p
  \end{bmatrix}
  =
  \begin{bmatrix}
    Q_1R_1\\
    \vdots\\
    Q_pR_p
  \end{bmatrix}.
\end{align}
The word local emphasizes that the problems are so small that they can be solved
without communication.

Then, the QR factorization of the vertically stacked $R$-factors
is computed
\begin{align}
  \label{eq:reducedQR}
  \begin{bmatrix}
    R_1\\
    \vdots\\
    R_p
  \end{bmatrix} =
  \begin{bmatrix}
    S_1\\
    \vdots\\
    S_p
  \end{bmatrix}R
\end{align}
and the global QR factorization of $A$ is then given by
\begin{align}
  A =
  \begin{bmatrix}
    Q_1S_1\\
    \vdots\\
    Q_pS_p
  \end{bmatrix}R.
\end{align}
To avoid that the QR factorization of the $R$-factors grows to large and thus
become too costly, this
algorithm is applied recursively, meaning it is applied itself to compute
the local QR factorizations in equation \eqref{eq:localqr}, or to compute the
reduced QR factorization \eqref{eq:reducedQR}.

Now we apply the TSQR algorithm to compute the factorization \eqref{eq:pqr}.
For that we assume that the matrix $Q$ was also computed with this algorithm and
is therefore stored in the TSQR representation, i.e.
\begin{align}
  Q = \begin{bmatrix}
    Q_1\\
    & Q_2\\
    && \ddots\\
    &&& Q_p
  \end{bmatrix}\begin{bmatrix}
    R_1\\
    R_2\\
    \vdots\\
    R_p
  \end{bmatrix},
\end{align}
where $Q_1,\ldots,Q_p$ are local orthogonal matrices and $R_1,\ldots,R_p$
are local upper triangular matrices.
As the matrix $Q$ is orthogonal we see that the stacked matrix of $R$-factors is
orthogonal too
\begin{align}
  \begin{bmatrix}
    R_1^T,\ldots,R_p^T
  \end{bmatrix}
  \begin{bmatrix}
    R_1\\
    \vdots\\
    R_p
  \end{bmatrix}
  & =
    \begin{bmatrix}
      R_1^T,\ldots,R_p^T
    \end{bmatrix}
    \begin{bmatrix}
      Q_1^T\\
      & \ddots\\
      && Q_p^T
    \end{bmatrix}
          \begin{bmatrix}
            Q_1\\
            & \ddots\\
            && Q_p
          \end{bmatrix}
                \begin{bmatrix}
                  R_1\\
                  \vdots\\
                  R_p
                \end{bmatrix}\\
  & = Q^TQ \\
  &= I.
\end{align}

Applying the TSQR algorithm on problem \eqref{eq:pqr} then means that we first
must solve the local QR factorizations
\begin{align}
  \begin{bmatrix}
    Q_i R_i & X_i
  \end{bmatrix}
  & =
    \begin{bmatrix}
      Q_i& U_i
    \end{bmatrix}
           \begin{bmatrix}
             R_i & P_i\\
             0 & N_i
           \end{bmatrix}
\end{align}
for every $i=1\ldots,p$,
which is equivalent to
\begin{align}
    \begin{bmatrix}
    Q_i & X_i
  \end{bmatrix}
  & =
    \begin{bmatrix}
      Q_i& U_i
    \end{bmatrix}
           \begin{bmatrix}
             I & P_i\\
             0 & N_i
           \end{bmatrix}.
\end{align}
This problem is the same as \eqref{eq:pqr} but of smaller size and can be
solved locally by one of the previous mentioned methods (or by a TSPQR algorithm
itself to introduce recursion).
We call the algorithm that is used for solving these problems the local subalgorithm.

Once all the local problems are solved, a \emph{reduction step} computes the QR
factorization of the stacked $R$ factors:
\begin{align}
  \begin{bmatrix}
    R_1 & P_1\\
    0 & N_1\\
    \vdots & \vdots\\
    R_p & P_p\\
    0 & N_p\\
  \end{bmatrix}
        &= \begin{bmatrix}
          R_1 & \widetilde{P}_1\\
          0 & \widetilde{N}_1\\
          \vdots\\
          R_p & \widetilde{P}_p\\
          0 & \widetilde{N}_p\\
        \end{bmatrix}
  \begin{bmatrix}
    I & P\\
    0 & N
  \end{bmatrix},
\end{align}
which is again a variant of problem \eqref{eq:pqr}, as we have seen that the
stacked R-factors are orthogonal.
The algorithm that is used to solve this problem is called the reduction subalgorithm.

The solution of the global system is then given by
\begin{align}
  \begin{bmatrix}
    Q & X
  \end{bmatrix}
  &= \begin{bmatrix}
      \begin{bmatrix}
        Q_1 & U_1
      \end{bmatrix}
      \begin{bmatrix}
        R_1 & \widetilde{P}_1\\
        0 & \widetilde{N}_1
      \end{bmatrix}\\
      \vdots\\
      \begin{bmatrix}
        Q_p & U_p
      \end{bmatrix}
      \begin{bmatrix}
        R_p & \widetilde{P}_p\\
        0 & \widetilde{N}_p
      \end{bmatrix}\\
    \end{bmatrix}
  \begin{bmatrix}
    I & P\\
    0 & N
  \end{bmatrix}.
\end{align}
The matrix $U$ can then be computed by block-wise matrix products
\begin{align}
  U = \begin{bmatrix}
    \begin{bmatrix}
        Q_1 & U_1
      \end{bmatrix}
      \begin{bmatrix}
        \widetilde{P}_1\\
        \widetilde{N}_1
      \end{bmatrix}\\
      \vdots\\
      \begin{bmatrix}
        Q_p & U_p
      \end{bmatrix}
      \begin{bmatrix}
        \widetilde{P}_p\\
        \widetilde{N}_p
      \end{bmatrix}
    \end{bmatrix}.
\end{align}

\begin{algorithm}
  \caption{\methodname{TreeTSPQR}}
  \label{alg:treetspqr}
  \begin{algorithmic}
    \For{$i = 1,\ldots,p$}
    \State Solve $\begin{bmatrix}
      Q_i & X_i
    \end{bmatrix} =
    \begin{bmatrix}
      Q_i & U_i
    \end{bmatrix}
    \begin{bmatrix}
      I & P_i\\
      & N_i
    \end{bmatrix}$
    \EndFor
    \State Solve $
    \begin{bmatrix}
      R_1 & P_1\\
      & N_1\\
      \vdots & \vdots\\
      R_p & P_p\\
      & N_p
    \end{bmatrix} =
    \begin{bmatrix}
      R_1 & \widetilde{P}_1\\
      & \widetilde{N}_1\\
      \vdots & \vdots\\
      R_p & \widetilde{P}_p\\
      & \widetilde{N}_p\\
    \end{bmatrix}
    \begin{bmatrix}
      R & P\\
        & N
    \end{bmatrix}$
  \end{algorithmic}
\end{algorithm}

\cref{alg:treetspqr} shows the pseudocode for the tree variant of the TSPQR
algorithm (\methodname{TreeTSPQR}).
As the local subproblems are independent they can be solved in parallel.
The size of the subproblems should be chosen such that all data fits in the
cache, so it must only be loaded once from the main memory.
In \cref{sec:numeric} (\cref{fig:local_size_benchmark}) we show a benchmark regarding the size of the subproblems.

\subsection{FlatTSPQR}
There is also a flat variant of the TSQR algorithm that we apply on problem
\eqref{eq:pqr} as well, to derive a flat variant of the TSPQR algorithm.
Flat TSQR computes a QR decomposition by factorizing the upper part of the
matrix and proceed with the remainder stacked with the R-factor of the upper part
\begin{equation}
  \begin{aligned}
  A &=
  \begin{bmatrix}
     Q_1R_1\\
     A_1
  \end{bmatrix} =
  \begin{bmatrix}
    Q_1\\
    &I
  \end{bmatrix}
      \begin{bmatrix}
        R_1\\
        A_1
      \end{bmatrix}
      =
      \begin{bmatrix}
        Q_1\\
        &I
      \end{bmatrix}
          \begin{bmatrix}
            Q_2\\
            &I
          \end{bmatrix}
              \begin{bmatrix}
        R_2\\
        A_2
      \end{bmatrix} = \ldots
  \\
  \label{eq:flat_tsqr}
  &=
    \begin{bmatrix}
      Q_1\\
      &I
    \end{bmatrix}
        \cdots
        \begin{bmatrix}
          Q_{p-1}\\
          &I
        \end{bmatrix}
            Q_p R.
          \end{aligned}
        \end{equation}
        In this representation all factors $
\begin{bmatrix}
  Q_i\\
  & I
\end{bmatrix}$ varying in size.

Applying this method to problem \eqref{eq:pqr} leading to the following method.
We assume that $Q$ is stored as a product of the form
\begin{align}
  Q = \begin{bmatrix}
      Q_1\\
      &I
    \end{bmatrix}
        \cdots
        \begin{bmatrix}
          Q_{p-1}\\
          &I
        \end{bmatrix}
            Q_p
\end{align}
and like in equation \eqref{eq:flat_tsqr} the
matrix $X$ is partitioned accordingly
\begin{align}
  X =
  \begin{bmatrix}
    X_1\\
    \vdots\\
    X_p
  \end{bmatrix}.
\end{align}
Starting with solving $\begin{bmatrix}Q_1& X_1
\end{bmatrix}
= \begin{bmatrix}Q_1& U_1
\end{bmatrix}\begin{bmatrix}I& P_1\\
  0 & N_1
\end{bmatrix}$, this leads to the sequence of problems
\begin{align}
  \begin{bmatrix}\multirow{3}{*}{$Q_i$}& P_{i-1}\\
    & N_{i-1}\\
    & X_i
  \end{bmatrix}
      = \begin{bmatrix}Q_i& U_i
      \end{bmatrix}\begin{bmatrix}I& P_i\\
        0 & N_i
      \end{bmatrix}
\end{align}
for $i=2,\ldots,p$.
Once all these problems are solved, \eqref{eq:pqr}
can be written as
\begin{align}
  \label{eq:flattspqr_q}
  \begin{bmatrix}
    Q & X
  \end{bmatrix}
  &=
    \begin{bmatrix}
      Q_1 & U_1\\
      && I
    \end{bmatrix}
         \cdots
         \begin{bmatrix}
           Q_{p-1} & U_{p-1}\\
           && I
         \end{bmatrix}
         \begin{bmatrix}
           Q_p & U_p\\
         \end{bmatrix}
  \begin{bmatrix}
    I & P_p\\
      & N_p
  \end{bmatrix}.
\end{align}
Here the matrix $U$ is assembled by computing the following matrix products
\begin{align}
  U = \begin{bmatrix}
      Q_1 & U_1\\
      && I
    \end{bmatrix}
         \cdots
         \begin{bmatrix}
           Q_{p-1} & U_{p-1}\\
           && I
         \end{bmatrix}
           U_p.
\end{align}

\begin{algorithm}
  \caption{\methodname{FlatTSPQR}}
  \label{alg:flattspqr}
  \begin{algorithmic}
    \State Solve $\begin{bmatrix}
      Q_1 & X_1
    \end{bmatrix} =
    \begin{bmatrix}
      Q_1 & U_1
    \end{bmatrix}
    \begin{bmatrix}
      I & P_1\\
      & N_1
    \end{bmatrix}$
    \For{$i=2,\ldots,p$}
    \State Solve $\begin{bmatrix}
      \multirow{3}{*}{$Q_i$} & P_{i-1}\\
      & N_{i-1}\\
      & X_i
    \end{bmatrix} =
    \begin{bmatrix}
      Q_i & U_i
    \end{bmatrix}
    \begin{bmatrix}
      I & P_i\\
      & N_i
    \end{bmatrix}$
    \EndFor
    \State $P = P_p\qquad N = N_p$
  \end{algorithmic}
\end{algorithm}
\cref{alg:flattspqr} shows the \methodname{FlatTSPQR} algorithm.
In contrast to the \methodname{TreeTSPQR} algorithm the local subproblems are
not independent and hence cannot be solved in parallel.

Like the TSQR algorithm, both variants inherit the stability properties of the
used subalgorithms. In the tree variant both choices of subalgorithms,
local and reduction, must be stable to
obtain a stable method.
We verify this numerically in \cref{sec:numeric} (\cref{fig:heatmap}).
Compared with the \methodname{BCGS-PIP} algorithm all algorithms load the
matrix $Q$ and $X$ only twice from the main memory but the TSPQR algorithms
are stable (if a proper subalgorithm is used).
Furthermore, both algorithms, \methodname{BCGS-PIP} and \methodname{TreeTSPQR}, only used one
synchronization point in a parallel setting.
However, we will see in \cref{sec:framework} that it is faster
to use \methodname{BCGS-PIP+} for the reduction in the message passing level, instead of
using a reduction tree.
\section{TSPQR as a Framework}
\label{sec:framework}
In this section we illustrate how the TSPQR methods can be used as a
framework to combine the building blocks and design a method tailored on a
specific hardware.
In particular for the \methodname{TreeTSPQR} method different sub-methods can be chosen for
computing the local and reduced problems.
As an example we use a CPU based cluster, consisting of multiple nodes with
multiple cores organized with a cache hierarchy.
The nodes are connected by a network.
Other architectures like GPUs or more sophisticated network topologies can be
treated similarly.

To investigate the choice of methods and assess the performance we employ
performance models, presented in the following subsections.
For our example these are the roofline model \cite{williams2008roofline} to model the intra-node
communication and the LogP model~\cite{culler1993logp} to model message passing performance.

\subsection{Intra-Node (Roofline Model)}
\label{sec:inter_node_performance}
To estimate the runtime of an algorithm on a single processor, we use the
roofline model \cite{williams2008roofline}.
It models the hardware by the peak performance $\pi$ [$\si{\flop\per\second}$]
and its peak memory bandwidth $\beta$ [$\si{\byte\per\second}$].
An algorithm specifies the two parameters amount of data $\delta$ [$\si{\byte}$]
and computation effort $\gamma$ [$\si{\flop}$].
The estimated runtime of the algorithm is then given by
\begin{align}
  \rho = \min\left(\frac{\gamma}{\pi}, \frac{\delta}{\beta}\right).
\end{align}
If the minimum is attained by the first term, the algorithm is
called \emph{compute-bound}, otherwise it is called \emph{memory-bound}.
As most modern architectures have separated channels for reading and writing to
the memory, we count only the data that is read.
In the following we determine the parameters $\gamma$ and $\delta$ for the
derived algorithms.
We only perform this analysis for stable methods, namely \methodname{BCGS-PIP+},
\methodname{HH} and TSPQR methods build up on them.

We consider two stages of the algorithms.
The first stage extends the basis $Q$, such that the new Krylov dimensions are
contained.
In the second stage, the matrix $U$ is assembled or a product $QC$ with matrix
$C\in\mathbb{R}^{(k+s)\times m}$ is computed.
This differentiation is useful for the analysis of the TSPQR performance.

\subsubsection{BCGS-PIP}
The \methodname{BCGS-PIP} algorithm performs $2n(sk+s^2)\,\si{\flop}$ computing the Gram
matrices $G$ and $P$, during that it reads the matrices $Q$ and $X$ from the memory.
We neglect terms for computing and factorizing the small matrices, as they do
not depend on $n$.
For assembling $U$ the algorithm reads again $Q$ and $X$ from the memory and
performs $2ns^2+2n(sk+s^2)\,\si{\flop}$.
In total $\gamma_1 = (4sk+6s^2)n\,\si{\flop}$ are executed and $Q$ and $X$ are read
two times from the memory each. Hence $\beta_1 = 16(k+s)n\,\si{\byte}$ are read from the
memory using double precision numbers ($8\,\si{byte}$).

The \methodname{BCGS-PIP+} algorithm repeats the algorithm, hence
the \methodname{BCGS-PIP+} algorithm executes $\gamma_1 = (8sk+12s^2)n\,\si{\flop}$ and read
$\beta_1 = 32(k+s)n\,\si{\byte}$.

As the basis is stored explicitly as a orthogonal matrix, the assembly of $U$
comes for free.
However, if the method is used for local orthogonalization in a TSPQR algorithm, the matrix product
$\begin{bmatrix}Q_i&U_i
\end{bmatrix}C$ must be computed to assemble $U$, which performs $\gamma_2 =
2(k+s)mn\,\si{\flop}$ and reads again $\beta_2 = 8(k+s)n\,\si{\byte}$.

\subsubsection{Householder}
The main operation in the Householder method is to apply Householder reflectors.
Applying a Householder reflector $H$ \eqref{eq:hh} on a $n\times s$ matrix $X$ performs
$4ns\,\si{\flop}$ and loads the vector $v$ as well as $X$ two times, that are
$16(s+1)n\,\si{\byte}$.

The \methodname{HH} algorithm applies $k$ Householder
reflectors on the vector $X$ at the beginning.
This performs $4ksn\,\si{\flop}$ and transfers $16k(s+1)n\,\si{\byte}$.
In the $i$th loop iteration then $i$ Householder reflectors are applied on the
$i$th column of $X$.
In total this are $\sum_{i=0}^{s-1} 4in = 4\frac{s(s-1)}{2}n = 2s^2n -
2sn\,\si{\flop}$ and $\sum_{i=0}^{s-1}16i(1+1)n = 16s^2n - 16sn$ bytes are transferred.
Hence for the first phase $\gamma_1 = (4ks + 2s^2 - 2s)n\,\si{\flop}$ are
performed and $\beta_1 = 16(k(s + 1) + s^2 - s)n\,\si{\byte}$ are transferred.

To assemble the result $U$, all $k+s$ Householder reflectors are applied on the
first $s$ unit vectors.
This performs $\gamma_2 = 4(k+s)sn\,\si{\flop}$ and transfers $\beta_2 = 16(k+s)(s+1)n\,\si{\byte}$.

We see the great disadvantage of the Householder method, i.e. that the memory
transfers scale are of order $\mathcal{O}(ks + s^2)$ instead of $\mathcal{O}(k+s)$.

\subsubsection{TSPQR}
The principle of the TSPQR algorithm is to choose the local problem size so
small such that the local problem can be solved in cache.
Therefore the matrices $Q$ and $X$ are only loaded once from the memory for
computing the factorization in the first stage.
In the Arnoldi process then the matrix $U$ must be assembled explicitly as the
second stage. For that the matrices $Q$ and $U$ are read again to compute the
product with the first unit vectors.
In total the algorithm read each matrix twice.
Hence, as the \methodname{BCGS-PIP} algorithm, it reads $\beta_1 = \beta_2 =
8(k+s)n\,\si{\byte}$ in every stage.

Assuming that the algorithm uses local problems of size
$\bar{n}$ and the effort for solving the first stage of the local problem is
$\nu \bar{n}\,\si{\flop}$, the computational costs for solving all first
stages of all local problems in the \methodname{TreeTSPQR} method is $p \nu \bar{n}\,\si{\flop}$.
As we have $p\bar{n} = n$, the total effort is $\gamma_1 = \nu n\,\si{\flop}$.

Analogous we assume that it needs $\theta\bar{n}\,\si{\flop}$ to proceed the matrix
multiplication
\begin{align}
  \begin{bmatrix}Q_i&U_i
  \end{bmatrix}
                      \begin{bmatrix}
  P_i\\N_i
\end{bmatrix}.
\end{align}
Hence the effort to assemble $U$ is $\gamma_2 = p \theta \bar{n} = \theta n\,\si{\flop}$.
The performance computations for the \methodname{FlatTSPQR} can be made analogously.

In case of the \methodname{BCGS-PIP+} method used for the local orthogonalization, we have
$\nu = 8sk + 12s^2$ and $\theta = 2(sk + s^2)$.
In total it is $10sk + 14s^2\,\si{\flop}$.
In contrast, if we use \methodname{HH} method for local orthogonalization these
value are given by
$\nu = 4sk + 2s^2$ and $\theta = 4(sk+s^2)$, leading to $8sk + 6s^2\,\si{\flop}$.

Compared to the stable \methodname{BCGS-PIP+} algorithm, we expect that the TSPQR algorithms
using the \methodname{HH} method for local orthogonalization are faster,
as they read less data from the memory and executes less floating point operations.
This is confirmed in the numerical test in \cref{sec:numeric}.

\subsection{Inter-Node (Message Passing)}
\label{sec:message_passing_model}
To model the performance of the inter node communication, we follow the
LogP model \cite{culler1993logp}. It assumes that the time
to send a message with payload [$\si{\byte}$] of $d$ is
\begin{align}
  \frac{d}{\beta} + \alpha,
\end{align}
where $\beta$ is the bandwidth [$\si{\byte\per\second}$] and $\alpha$ is the
latency [$\si{\second}$] of the network.

In a all-reduce operation the nodes are organized in a tree where messages are
sent from the leafs to the root.
On every node in the tree a reduction operation is performed, reducing the
incoming messages. The result is then send to the parent node.
Following the LogP model, we model the execution time of a reduction communication as
\begin{align}
  \log(P)\left(\omega + \frac{d}{\beta} + \alpha\right),
\end{align}
where $P$ is the number of nodes and $\omega$ is the time [$\si{\second}$] to perform the
reduction on one tree node.

In both cases, \methodname{BCGS-PIP} and TSPQR, the message is a $s\times(k+s)$ matrix, i.e.\
$d=8(ks+s^2)\,\si{\byte}$.
The difference between \methodname{BCGS-PIP+} and TSPQR is the reduction operation.
In case of the \methodname{BCGS-PIP+} this reduction is summation.
Therefore we have $\omega_{+} = \mathcal{O}\left(sk+s^2\right)$.
In TSPQR, the reduction operation is solving the project and normalize problem
for the stacked matrices, hence $\omega_{\operatorname{TSPQR}} =
\mathcal{O}\left(sk^2 + s^2k\right)$.
We see that the reduction operation of the TSPQR algorithm is more
expensive by a factor of $k$.
Furthermore, the summation can be implemented using \texttt{MPI\_Allreduce}
which is an optimized implementation by the MPI vendor. For the TSPQR such an
MPI function does not exists, as the tree-nodes need to store the state of the
reduced basis $Q_i$ between calls.
Hence we have build our own implementation which is not much tuned.

Whether \methodname{BCGS-PIP+} or \methodname{TreeTSPQR} is faster depends on the network parameters $\alpha$
and $\beta$ as well as on $s$ and $k$.
From the theoretical site, \methodname{BCGS-PIP+} is faster as long as
\begin{equation}
  \begin{aligned}
    &&  2\log(P)\left(\omega_+ + d\beta + \alpha\right) &<
    \log(P)\left(\omega_{\operatorname{TSPQR}} + d\beta + \alpha\right)\\
    \Leftrightarrow && 2\omega_{+} + d\beta + \alpha &<
    \omega_{\operatorname{TSPQR}}.
  \end{aligned}
\end{equation}
In our test setting that holds almost always true, as the interconnect is quite
fast.
However, in slow networks, e.g.\ if the latency is larger than the duration of
the reduction
operation, $\alpha \gg \omega_{\operatorname{TSPQR}}$,
\methodname{TreeTSPQR} can be a better choice.

\begin{figure}
  \centering
  \includegraphics[width=\textwidth]{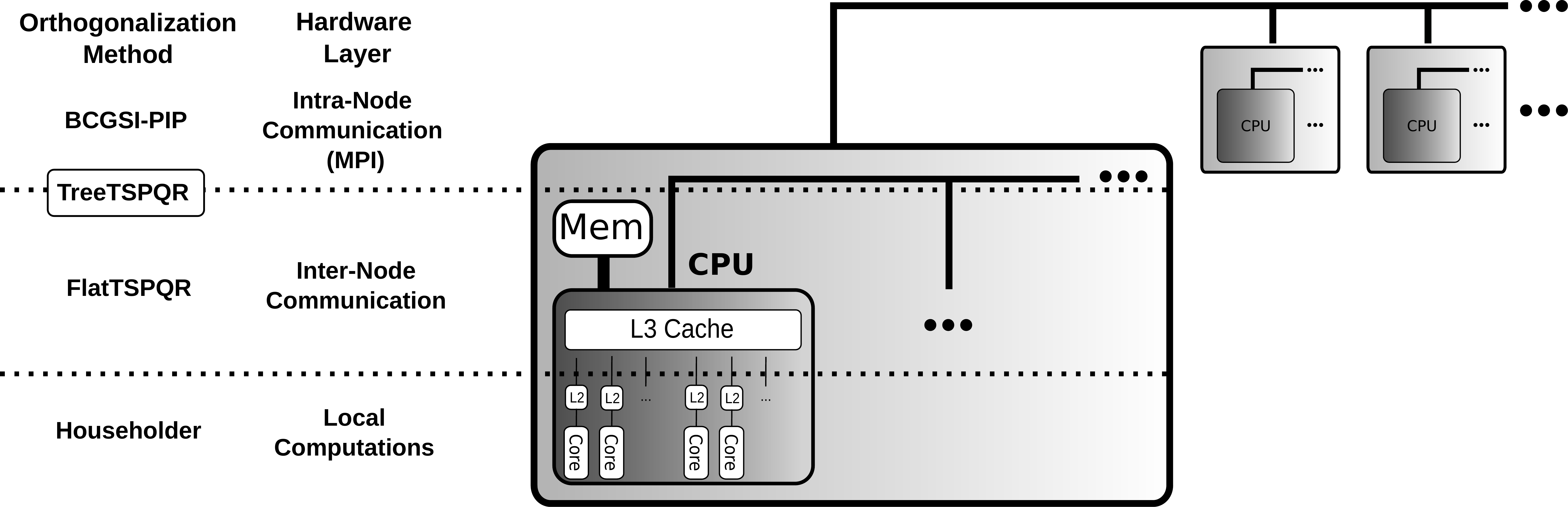}
  \caption{Design of the algorithm to the computer architecture. \methodname{TreeTSPQR}
    connects the \methodname{FlatTSPQR} algorithm for the computations on the nodes and the
    \methodname{BCGS-PIP+} algorithm that is used for the reduction on the message passing level.}
  \label{fig:architecture}
\end{figure}

That means that TSPQR performs better for the solution of the problem on one
node, while \methodname{BCGS-PIP+} performs better for message passing environments.
We therefore propose to combine the methods and use \methodname{FlatTSPQR} using
\methodname{HH} for the local orthogonalization to solve the
problem on one MPI rank and do the reduction of the results using the \methodname{TreeTSPQR}
method, using the \methodname{BCGS-PIP+} method for the reduction problem.
This enables to utilize both advantages - good node performance and using the
optimized inter-node communication pattern of the \methodname{BCGS-PIP+} method.
The design of the method and its adaption to the hardware is illustrated in
\cref{fig:architecture}.

\section{Numerical Experiments}
\label{sec:numeric}
For the benchmarks and stability experiments we implemented the methods using
the Eigen \cpp framework~\cite{eigen}.
We used MPI for the parallelization using one rank per core.
The source code is provided as supplementary material.

\subsection{Stability}
To investigate the stability of the algorithms, we use the Stewart matrices
presented in \cite{stewart2008block} to construct matrices with a given condition number.
These matrices are constructed from a random matrix $A$, that is singular value
decomposed $A=U\Sigma V^T$.
The orthogonal matrices from the singular value decomposition are then recombined
to the desired matrix $X = U\tilde{\Sigma}V^T$, where $\tilde{\Sigma}$ is a
diagonal matrix with exponentially increasing diagonal entries $\frac{1}{\kappa},\ldots,1$, such that $X$ has
the desired condition number $\kappa$.
If not otherwise stated all experiments are carried out with $n=2^{16}, {k = 32}$ and
$s=4$. For the TSPQR methods the problem is subdivided into problems of size
$256$ rows.

\begin{figure}
  \centering
  \subfloat[basic variants]{
    \includegraphics[width=0.47\textwidth]{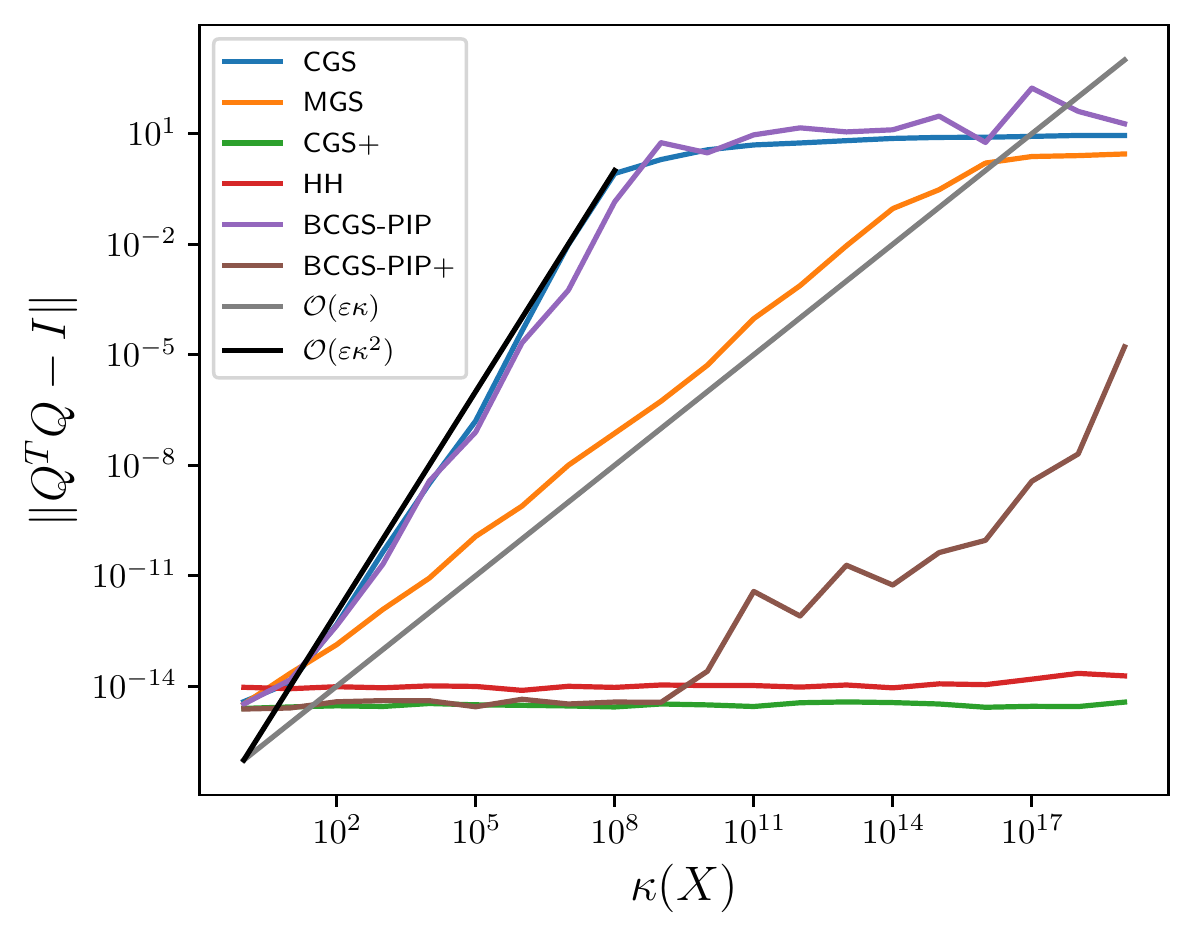}
    \label{fig:kappa_a}
  }
  \subfloat[TSPQR]{
    \includegraphics[width=0.47\textwidth]{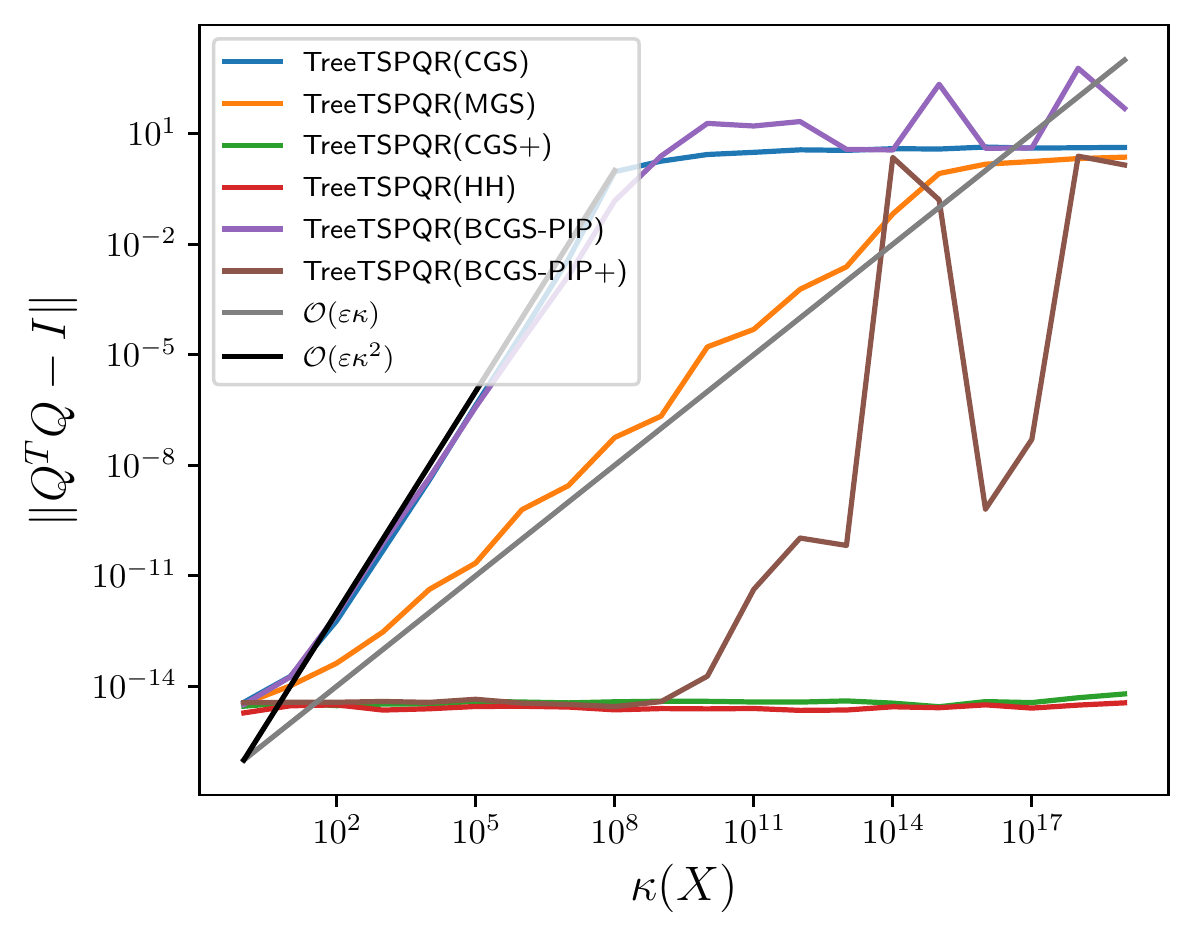}
    \label{fig:kappa_b}
  }
  \caption{Orthogonalization errors of presented algorithms for Stewart matrices
    with different condition numbers. Colors encode the different methods. The
    black and gray line mark $\mathcal{O}(\varepsilon\kappa)$
    and $\mathcal{O}(\varepsilon\kappa^2)$.}
  \label{fig:kappa}
\end{figure}
The first experiment is inspired by Carson et al.~\cite{carson2020overview}.
For that we compute the QR factorization of a Stewart matrix with different
condition numbers using the presented algorithms applied on the $s$-block
columns (\cref{alg:qr}).
\cref{fig:kappa_a} shows the orthogonalization error $\varepsilon_\perp$ as
defined in equation \eqref{eq:ortho_error} for the resulting $Q$-factor.
We see that the traditional methods, Gram-Schmidt and Householder, behave as
expected.
The orthogonalization errors of the reiterated classical Gram-Schmidt
method (\methodname{CGS+}) and the \methodname{HH} method are of order $\mathcal{O}(\varepsilon)$.
For the modified Gram-Schmidt method the orthogonalization error is of order
$\mathcal{O}(\kappa\varepsilon)$ and for the classical Gram-Schmidt method
\methodname{CGS} it is of
order $\mathcal{O}(\kappa^2\varepsilon)$.

The orthogonalization error of the \methodname{BCGS-PIP} method is of order
$\mathcal{O}(\kappa^2\varepsilon)$, while its reiterated variant is stable up to
a condition number satisfying
\begin{align}
  \label{eq:BCGS_condition}
  \varepsilon\kappa^2 \leq \frac12.
\end{align}

The orthogonalization errors for different choices of subalgorithms algorithms
in the \methodname{TreeTSPQR} algorithm are plotted in \cref{fig:kappa_b}.
We see that the stability of the \methodname{TreeTSPQR} algorithm equals
the stability of the used subalgorithm.
The mentioned subalgorithm is used for the local as well as for the reduction
step.
The \methodname{TreeTSPQR} seems to be more sensible for condition
\eqref{eq:BCGS_condition} if \methodname{BCGS+} is used as the subalgorithm.

\begin{figure}
  \centering
  \includegraphics[width=.6\textwidth]{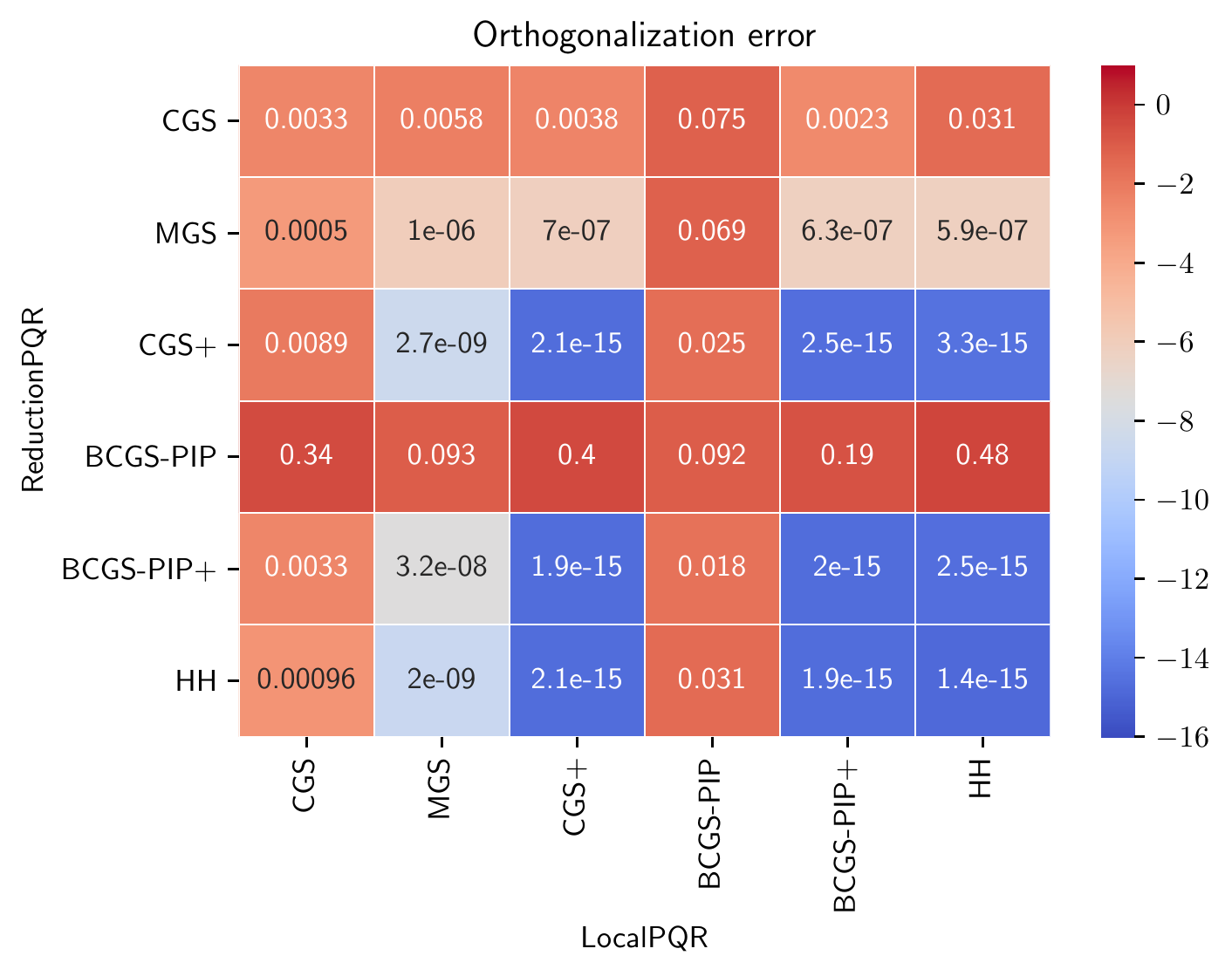}
  \caption{Orthogonalization errors for different combinations of local and
    reduction submethods used in \methodname{TreeTSPQR} for a Stewart matrix
    with condition number of $\kappa = \num{e8}$.}
  \label{fig:heatmap}
\end{figure}

In the next experiment we investigate how the orthogonalization error depends on
the choice of the local and reduction subalgorithm respectively.
For that we compute the QR factorization of a Stewart matrix with condition
number $\kappa = \num{e8}$ with different combinations of subalgorithms.
The result is shown in \cref{fig:heatmap}.
Red color indicates a high orthogonalization error, while blue color indicates
stability.
Columns encode the method that is used to solve the local problems and rows
encode the reduction orthogonalization method.
It shows that both subalgorithms must be stable to obtain a stable method.

\begin{table}
  \centering
  \footnotesize
  \setlength\tabcolsep{5pt} 
  \def\arraystretch{1.2} 
  \caption{Stability results for the \methodname{TreeTSPQR} algorithm with
    different number of recursion levels and different number of local problems
    (local problem size).}
  \subfloat[recursion levels]{
    \label{tab:recursion_depth_error}
    \centering
    \begin{tabular}{|c*{6}{|S[parse-numbers = true, table-format=1.1e-2,
      round-mode = places, round-precision = 1, exponent-product = \cdot,
      tight-spacing=true]}|}
      \hline
      & \multicolumn{3}{c|}{$e_\perp$} & \multicolumn{3}{c|}{$\|A-QR\|$}\\
      \hline
      level & \methodname{BCGS-PIP} & \methodname{BCGS-PIP+} & \methodname{Householder}& \methodname{BCGS-PIP} & \methodname{BCGS-PIP+} & \methodname{Householder}\\
      \hline\hline
      0&0.99749e-09&2.32836e-15&9.07492e-15&5.8348e-16& 6.70002e-16&4.07612e-15\\
      1&2.35387e-09&2.39143e-15&2.15442e-15&5.54934e-16&1.11875e-15&1.19797e-15\\
      2&7.58058e-09&3.60049e-15&1.94408e-15&7.16932e-16&9.23065e-16&1.2159e-15\\
      3&2.60133e-09&4.38641e-15&2.11658e-15&1.06734e-15&1.14409e-15&1.45764e-15\\
      4&1.46254e-09&3.1776e-15&1.69669e-15&9.23807e-16&1.67634e-15&1.47585e-15\\
      5&5.25672e-09&4.02172e-15&1.95351e-15&1.37734e-15&1.09762e-15&1.57498e-15\\
      6&3.84084e-09&4.35938e-15&1.98905e-15&1.23418e-15&1.42175e-15&1.82878e-15\\
      7&2.90005e-09&5.2614e-15&1.41694e-15&1.4143e-15& 2.27328e-15&1.95875e-15\\
      8&2.43657e-09&5.30096e-15&1.90436e-15&1.94375e-15&2.21207e-15&2.38217e-15\\
      \hline
    \end{tabular}
  }\\
  \subfloat[number of subproblems]{
    \centering
    \label{tab:local_size_error}
    \begin{tabular}{|c*{6}{|S[parse-numbers = true, table-format=1.1e-2,
      round-mode = places, round-precision = 1, exponent-product = \cdot,
      tight-spacing=true]}|}
            \hline
      & \multicolumn{3}{c|}{$e_\perp$} & \multicolumn{3}{c|}{$\|A-QR\|$}\\
      \hline
      sub\newline problems & \methodname{BCGS-PIP} & \methodname{BCGS-PIP+} & \methodname{Householder}& \methodname{BCGS-PIP} & \methodname{BCGS-PIP+} & \methodname{Householder}\\
      \hline\hline
      8   &3.03684e-09&1.76856e-15&2.07594e-15&6.49862e-16&7.14792e-16&9.89664e-16\\
      16  &3.33663e-09&2.58334e-15&1.69639e-15&8.59705e-16&6.43550e-16&8.27309e-16\\
      32  &3.67168e-09&2.63008e-15&1.81935e-15&6.06690e-16&7.27582e-16&8.16933e-16\\
      64  &4.74107e-09&2.25693e-15&2.19955e-15&6.11734e-16&7.96782e-16&1.42728e-15\\
      128 &5.9833e-09 &3.61907e-15&2.05285e-15&7.71861e-16&5.68321e-16&9.56823e-16\\
      256 &2.47485e-09&2.60008e-15&2.36896e-15&5.62788e-16&7.79341e-16&1.21995e-15\\
      512 &2.89839e-09&2.46391e-15&2.04322e-15&6.68429e-16&8.34059e-16&1.32417e-15\\
      1024&5.32721e-09&2.23119e-15&2.42144e-15&6.08644e-16&4.98339e-16&1.82894e-15\\
      \hline
    \end{tabular}
  }
\end{table}



In a further experiment, we take a look at the errors dependent on the numbers
of recursion levels in the \methodname{TreeTSPQR} algorithm and the number of processes used.
Here we use a condition number of $\kappa = \num{e4}$.
\Cref{tab:recursion_depth_error} shows the norm of the residual and the
orthogonalization error of the \methodname{TreeTSPQR} algorithm with different number of
recursion levels.
\Cref{tab:local_size_error} shows the residual norm and orthogonalization error
of the \methodname{TreeTSPQR} method with differently many subproblems in two levels.
The same system size was used, so if more subproblems are used the subproblems
are of smaller size.
We see that neither the orthogonality error nor the residual norm depend
crucially on the number of used recursion levels or number of subproblems.

\subsection{Performance}
We demonstrate the performance advantages of the presented orthogonalization
framework in this subsection.
First, we consider the performance on a single node where message passing is
cheap and the considerations in \cref{sec:inter_node_performance} are relevant.
These experiments are carried out on our AMD Epic 7501 compute server with 64
physical cores.
Each core has a $\SI{512}{\kilo\byte}$ L2 cache.
The cores are organized on 8 sockets with 8 cores each, where all cores on a
socket share an $\SI{8}{\mega\byte}$ L3 cache.

To compare only the performance of the orthogonalization procedure, we compare
the runtimes for computing a QR-factorization (\cref{alg:qr}) instead of the Arnoldi method.
With the Arnoldi method the timings would be unclear due to the application of
the operator.
We choose the problem size such that the L3 cache is exhausted.
The input matrix for the QR-decomposition has $64$ columns (fix) and
$2^{18}$ rows per process (weak scaling).
This leads to a problem size of $64\cdot 2^{18}P \cdot 8\si{\byte} =
P\cdot\SI{128}{\mega\byte}$.
We used a chunk size of $s=4$ and TSPQR methods use a local problem size of $n_p
= 2^{13}$ rows ($256\,\si{\kilo\byte}$).

\begin{figure}
  \centering
  \subfloat[Sequential]{
    \includegraphics[width=.47\textwidth]{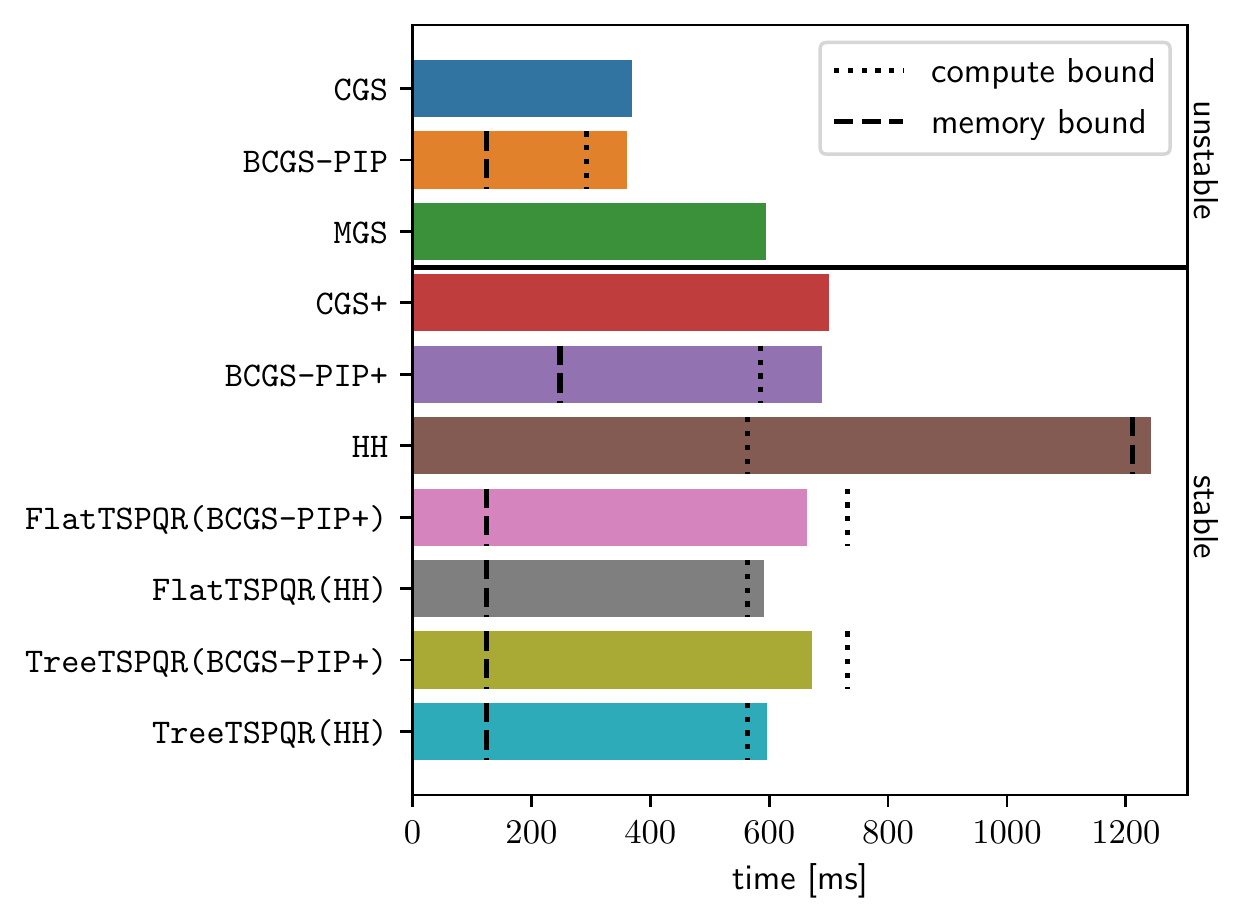}
  }
  \subfloat[Parallel ($64$ processes)]{
    \includegraphics[width=.47\textwidth]{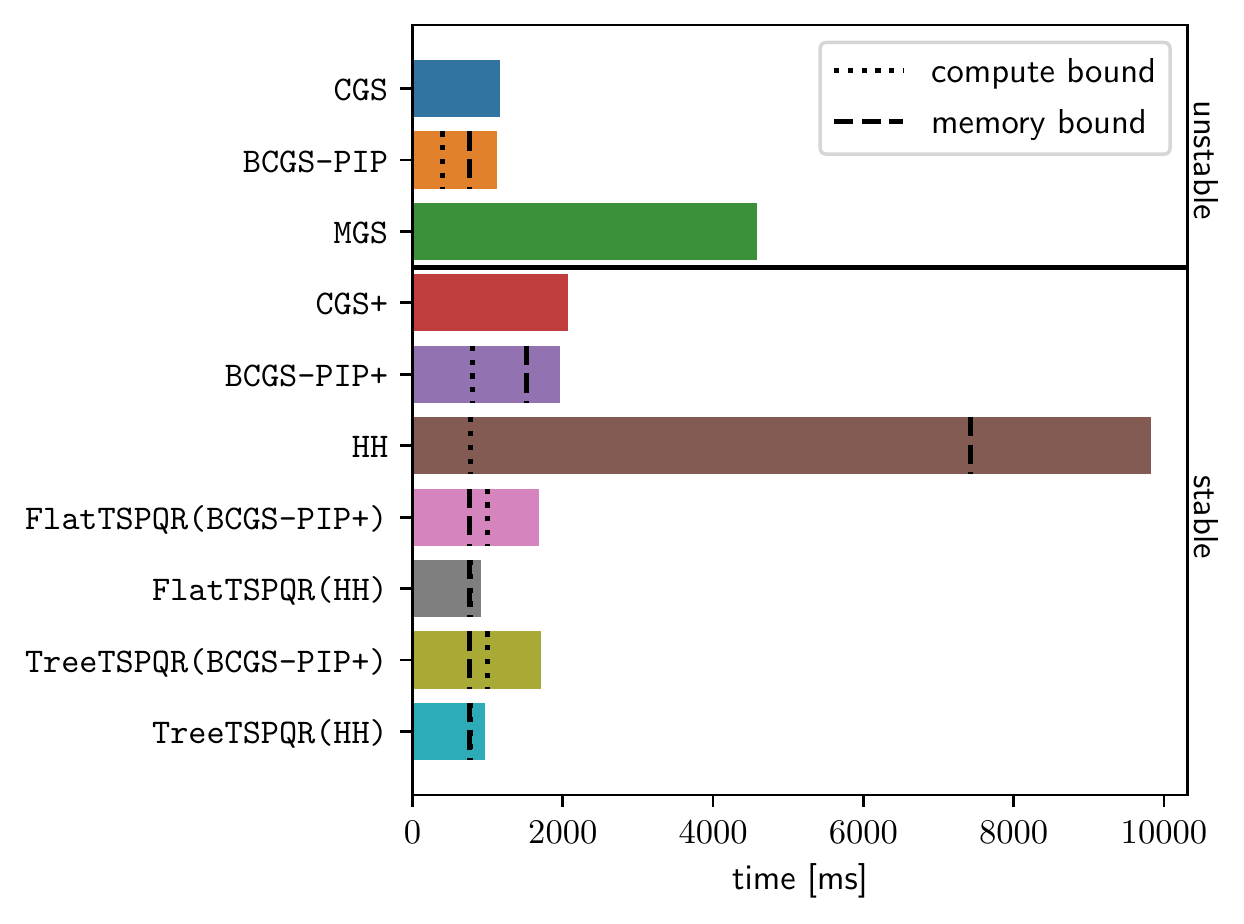}
  }
  \caption{Comparison of algorithm runtimes on a single node with on one core
    (sequential) and on all 64 cores (parallel). Dashed and dotted lines mark
    the prediction of the roofline model. The solid black line separates
    unstable from stable methods.}
  \label{fig:pqr_benchmark}
\end{figure}

In a first benchmark we compare the runtimes of the different
algorithms presented in this paper.
\Cref{fig:pqr_benchmark} shows the runtime of the different methods for a single
process and for $64$ processes.
At the top the unstable methods are displayed for comparison.
Black vertical lines mark the compute- (dotted) and memory-bound (dashed) that
is predicted by the roofline model in \cref{sec:inter_node_performance}.
To compute these bounds we measured the peak performance and memory bandwidth
using the likwid-bench tool \cite{likwid} and the \texttt{stream\_avx\_fma}
benchmark.
This benchmark is representative for the operations we perform in most
algorithms.
For one core we measured $\pi_1 = \SI{7.8}{\giga\flop\per\second}$ and $\beta_1
= \SI{18.4}{\giga\byte\per\second}$.
For the parallel test case on all $64$ cores we measured
$\pi_{64} = \SI{364}{\giga\flop\per\second}$ and
$\beta_{64} = \SI{192}{\giga\byte\per\second}$.

The prediction of the performance model matches the result.
The \methodname{HH} method is the slowest, due to the large amount of data that
needs to be transferred from the memory.
This issue is solved by the TSPQR methods which perform much better.
In both cases the TSPQR methods are the fastest stable methods.
In particular in the parallel case the TSPQR methods perform even faster as the
\methodname{BCGSI-PIP} method, while preserving the stability of the used
subalgorithm.
TSPQR methods that use \methodname{HH} as subalgorithm are faster as the one
used \methodname{BCGS-PIP+}. This is due to the higher computational effort of
the \methodname{BCGS-PIP+} method for the assembly of the output matrix $U$ if
it used in TSPQR methods, as discussed in \cref{sec:inter_node_performance}.

In the sequential case, the TSPQR methods using \methodname{BCGS-PIP+} as a
subalgorithm are even faster as predicted by the roofline model.
This is due to the fact, that this algorithm uses a lot of matrix-matrix products
which can be better optimized to achieve a higher flop-rate as the
\texttt{stream\_avx\_fma} benchmark.




In the next benchmarks, we investigate the best local problem size for the TSPQR
algorithms.
Naturally, we want to choose the number of rows of the local spaces such that the
local problem fits into the L2 cache.
Unfortunately, the local problem size grows when the Krylov space grows, so that
also the local problems grow.

\begin{figure}
  \centering
  \subfloat[Sequential]{
    \includegraphics[width=.47\textwidth]{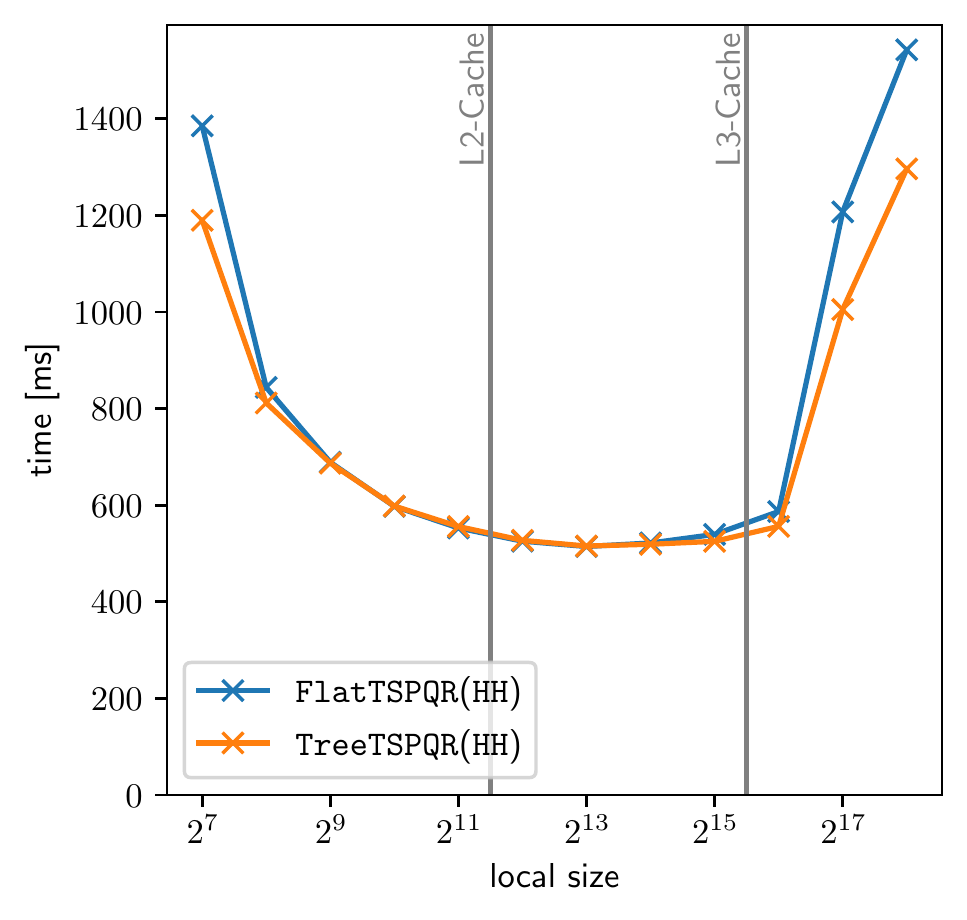}
  }
  \subfloat[Parallel ($64$ processes)]{
    \includegraphics[width=.47\textwidth]{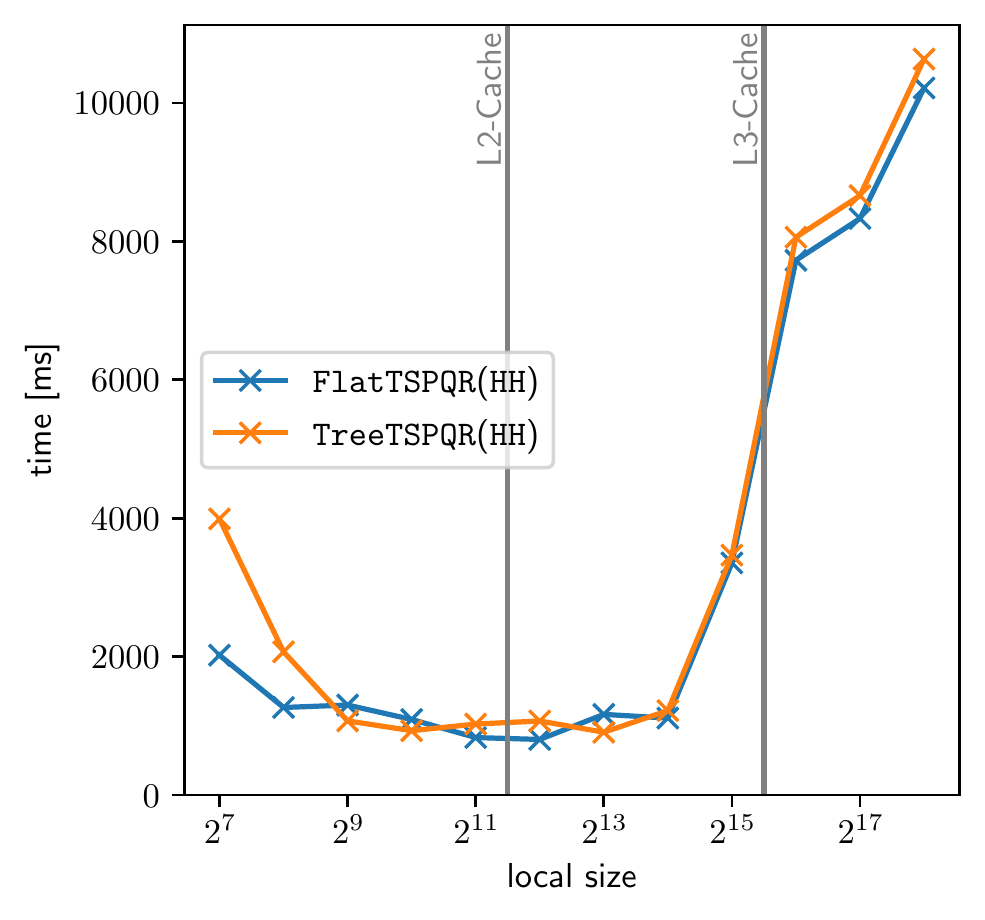}
  }
  \caption{Benchmark for different local problem size in TSPQR algorithms.
  Using $s=2^{18}$ rows, $k=64$ and $s=8$ columns on each process (weak
  scaling). Gray vertical lines mark the size of the L2- and L3- cache.}
  \label{fig:local_size_benchmark}
\end{figure}
In \cref{fig:local_size_benchmark} the runtimes for computing the QR
decomposition per number of local rows are
shown for the same setting as in \cref{fig:pqr_benchmark}.
The gray lines mark the size of the L2 and L3 cache for the largest case
($k=60, s=4$).
We see that the minimum runtime is approximately at $2^{12}$ in the sequential
and parallel case, which is slightly bigger than the L2 cache.
For fewer local rows the method introduces overhead that leads to higher
runtimes.
For more local rows, the local problem becomes bigger than the cache which then
leads to more communication between the memory and cache.
We see also that the local characteristic of the algorithm has far more impact
in the parallel setting.
This is due to the fact that multiple cores share memory bandwidth, which leads
to a smaller memory-bandwidth per process.

\begin{figure}
  \centering
  \resizebox{0.6\textwidth}{!}{
    \includegraphics[width=.8\textwidth]{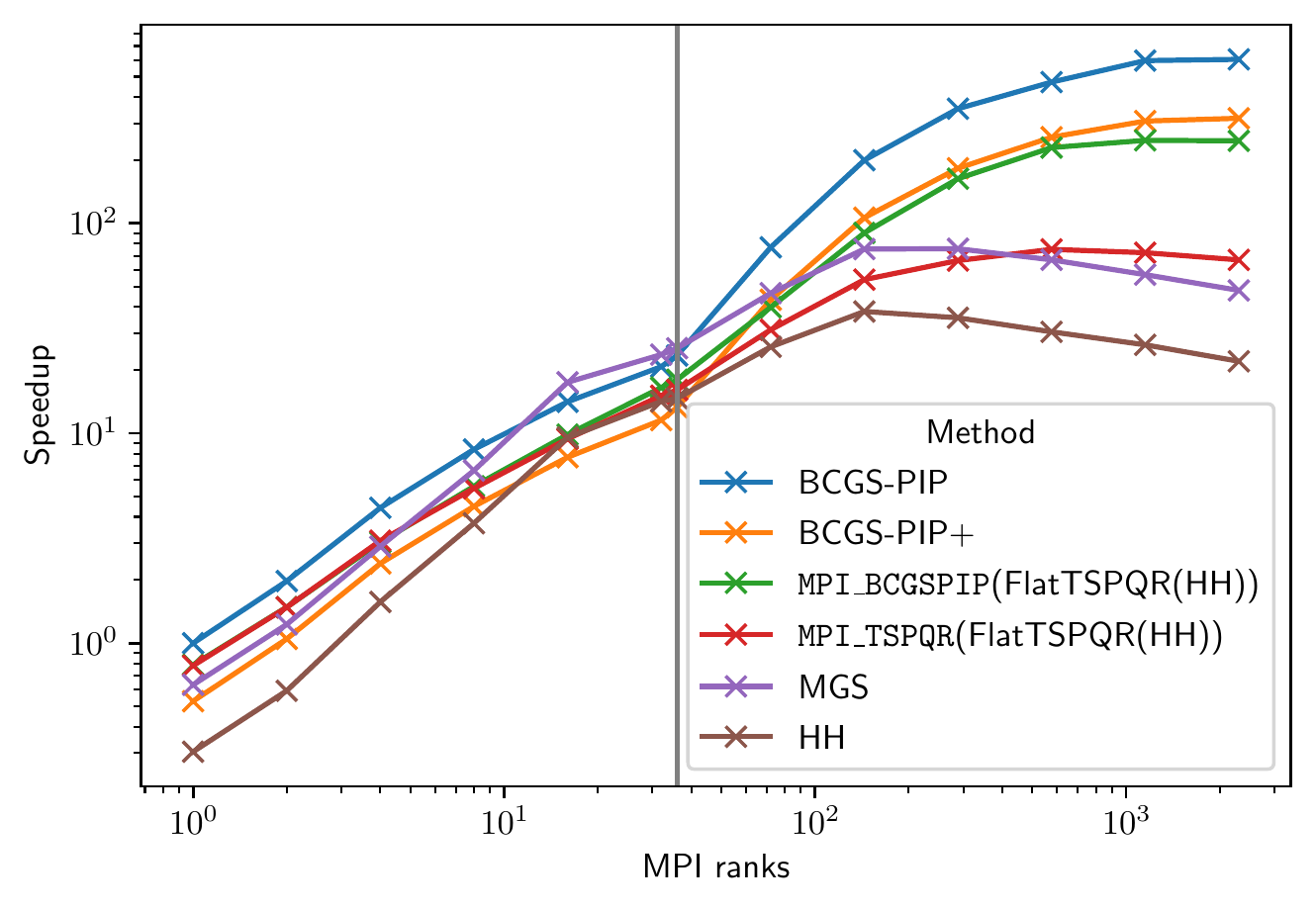}
  }
  \caption{Strong-scaling benchmark on the PALMAII super computer on up to 2304
    processes. The speedup is computed against the fastest method on one
    process (\methodname{BCGS-PIP}). Both axes are log-scale. Up to the gray
    vertical line the computation is on one node.}
  \label{fig:strong_scaling}
\end{figure}
Finally, we validated the performance expectations for the
message passing model from \cref{sec:message_passing_model}.
For that we used the supercomputer PALMAII of the university of Münster.
\cref{fig:strong_scaling} shows the speedup for the same problem as used in the
other benchmarks but in a strong scaling setting using up to 64 nodes with 36
processes each (2304 processes in total).
The vertical gray line marks 36 processes, which is the limit for which the
computation proceeds on one node.
At the scaling limit the problem size is $n = \frac{2^{18}}{2304}\approx 114$ per process.
We see that for large number of processes the performance of \methodname{BMGS}
and \methodname{HH} stagnate, as they need many synchronization points,
whereas the performance of
the \methodname{BCGS-PIP} methods is much better.
The runtime of the reiterated \methodname{BCGS-PIP+} method is almost twice the runtime of
the \methodname{BCGS-PIP} method.
In our framework, we implemented two variants for the reduction of the
\methodname{TreeTSPQR} method on the MPI
level, \methodname{BCGS-PIP+} and TSPQR, as discussed in \cref{sec:message_passing_model}.
In \cref{fig:strong_scaling} the labels are prefixed with \texttt{MPI\_} to
indicate the respective usage on the MPI level.
The \methodname{MPI\_TSPQR} implementation performs not as good as the \methodname{MPI\_BCGS-PIP+} methods as it
was predicted by out performance model.
But we see that the performance of the \methodname{MPI\_TSPQR} method is not
decaying as fast as for the \methodname{HH} or \methodname{MGS} algorithms.

\section{Conclusions and Outlook}
\label{sec:conclusions}
In this paper we presented a new orthogonalization framework TSPQR that
can be used to design orthogonalization algorithms for Krylov methods tailored
on a specific hardware.
The most common orthogonalization algorithms, Gram-Schmidt and Householder, are
used as building blocks to solve smaller, local problems, that can be solved
without communication.
In principle other algorithms can be used as well.
To demonstrate the usage and potential of the framework, we designed an
algorithm for a CPU-based HPC-Cluster combining the Householder algorithm used
for the local orthogonalization and the \methodname{BCGS-PIP+} algorithm for
reduction on the MPI layer.
Furthermore, we presented a performance analysis for the designed method and
presented numerical examples concerning the stability and performance.

The experiments showed that the novel framework can be used to design algorithms
that perform as good as the performance-optimal \methodname{BCGS-PIP} algorithm,
but preserve the stability of the used subalgorithm.
In addition, it provides more flexibility for the developers and enables them to
reuse existing implementations as building blocks.

The variety of orthogonalization algorithms that can be designed is manifold.
A future goal would be to try out other orthogonalization methods as building
blocks in this framework, e.g.\ paneled Householder algorithm
\cite{schreiber1989storage} and see whether the performance could be further
improved.
Another goal would be to design algorithms for other hardware like GPUs and
accelerators.
Furthermore, we showed the stability of the TSPQR algorithms only experimentally.
It would be worth to do an elaborate stability analysis to get further insights
into the framework.

For the future, the MPI standard specification
could be extended to allow stateful tree-reduction operations like it is needed by the
\methodname{TreeTSPQR} algorithm.
This would allow the vendor to optimize this operation on a specific hardware.
In a further step, it could be beneficial to introduce network hardware, that
implements the stateful tree-reduction where the state is stored on the network
switches to improve this operation.

\bibliographystyle{siamplain}
\bibliography{references}
\end{document}